# Some questions of connection between summation functions and the corresponding Dirichlet series

Victor Volfson

ABSTRACT. The paper proves a generalization of Wintner's theorem on the asymptotics of summation functions to the case of summation functions with nonlinear asymptotics. The class of arithmetic functions that have a logarithmic asymptotic mean is studied. The Kronecker lemma is generalized to the case when the corresponding Dirichlet series diverges. Several assertions are proven and illustrative examples are given.

## 1. INTRODUCTION

In general, an arithmetic function is a function defined on the set of natural numbers and taking values on the set of complex numbers. The name arithmetic function is due to the fact that this function expresses some arithmetic property of the natural series.

The summation function of an arithmetic function $f$ is called a function of the form:

$$M(f,x) = \sum_{n \leq x} f(n), \qquad (1.1)$$

where $x$ is a real number.

The Dirichlet series of an arithmetic function $f$ is called a function of the form:

$$F(s) = \sum_{n=1}^{\infty} \frac{f(n)}{n^s}, \qquad (1.2)$$

where $s$ is a complex number.

The Dirichlet convolution of arithmetic functions $f, g$ is the operation *, which is defined as follows:

$$(f * g)(n) = \sum_{d|n} f(d)g(n/d) = \sum_{ab=n} f(a)g(b). \qquad (1.3)$$

Wintner proved the theorem on the asymptotic mean value of an arithmetic function [1].

Suppose $f(n) = (1 * g)(n) = \sum_{d|n} g(d)$ (1 is an arithmetic function that takes the value 1 for all natural numbers $n$) and the condition is satisfied that the Dirichlet series of the arithmetic function $g$ is absolutely convergent, i.e. $\sum_{n=1}^{\infty} \frac{|g(n)|}{n} < \infty$.

Then the arithmetic function $f$ has an asymptotic mean value, which is equal to:

$$M(f) = \lim_{x \to \infty} \frac{M(f,x)}{x} = \lim_{x \to \infty} \frac{1}{x} \sum_{n \leq x} f(n) = \sum_{n=1}^{\infty} \frac{g(n)}{n}. \qquad (1.4)$$

A corollary of this theorem is that when these conditions are met, the following linear asymptotics is true for the summation function of the arithmetic function $f$:

$$M(f,x) = x\sum_{n=1}^{\infty} \frac{g(n)}{n} + o(x), x \to \infty. \qquad (1.5)$$

The linear asymptotic behavior of summation functions was studied in detail in [2].

It is interesting to generalize Wintner's theorem to the case of nonlinear asymptotics for the summation function of an arithmetic function $f$. This will be done in Chapter 2 of this work.

The nonlinear asymptotic behavior of summation functions is considered in detail in [3].

The assertion [4] is known. If $g'(\xi)$ has a constant sign ($g(\xi)$-monotone) and $g(x) \to 0$ at $x \to \infty$, then:

$$\sum_{a<n\leq x} g(n) = \int_a^x g(\xi)d\xi + C + O(|g(x)|), \qquad (1.6)$$

where $C = \int_a^{\infty} (\xi - [\xi])g'(\xi)d\xi$.

An example of using (1.6) is the asymptotics for the summation function:

$$\sum_{n\leq x} \frac{1}{n} = \ln x + \gamma + O(\frac{1}{x}),$$

where $\gamma$ is Euler's constant.

Recall that an arithmetic function $f$ has a logarithmic asymptotic mean at $x \to \infty$ equal to $A$ (constant) if:

$$\sum_{n\leq x} \frac{f(n)}{n} = A\ln x + O(\ln x), (x \to \infty). \qquad (1.7)$$

An arithmetic function $f$ has a strongly logarithmic asymptotic mean at $x \to \infty$ equal $A$ (constant) if:

$$\sum_{n \leq x} \frac{f(n)}{n} = A \ln x + B + o(1), (x \to \infty),  \qquad (1.8)$$

where $B$ is a constant.

Unfortunately, arithmetic functions are mostly not monotonic, so using formula (1.6) is often impossible. Therefore, it is of interest to find a class of arithmetic functions that have asymptotic behavior of the form (1.7) or (1.8). This class of arithmetic functions will be discussed in Chapter 3 of this work.

Kronecker's lemma is widely used in probability theory [5] and number theory [6]. Let there is an arithmetic function $f$ and a complex number $s$ ($\operatorname{Re} s > 0$) for which:

$$\sum_{n=1}^{\infty} \frac{f(n)}{n^s} < \infty, \qquad (1.9)$$

then $\lim_{x \to \infty} \dfrac{1}{x^s} \sum_{n \leq x} f(n) = 0$ \qquad (1.10)

A corollary of this lemma is that the following asymptotic behavior for the summation function is valid:

$$\sum_{n \leq x} f(n) = o(x^s), \qquad (1.11)$$

when these conditions are met.

Example. It is known that $\sum_{n=1}^{\infty} \dfrac{\sigma(n)}{n^s} = \zeta(s-1)\zeta(s) < \infty$ at $\operatorname{Re} s > 2$. Thus, when $s = 3$ we get $\sum_{n=1}^{\infty} \dfrac{\sigma(n)}{n^s} = \zeta(2)\zeta(3) < \infty$. Therefore, based on the corollary of Kronecker's lemma $\sum_{n \leq x} \sigma(n) = o(x^3)$.

However, Kronecker's lemma and its corollary do not answer the question of what the asymptotic behavior will be, for example, of the summation function $\sum_{n \leq x} \dfrac{\sigma(n)}{n^s}$ for $s = 1, 2$ when the corresponding Dirichlet series diverges.

A similar question arises for other summation functions of the form $\sum_{n \leq x} \frac{f(n)}{n^s}$, when the corresponding Dirichlet series diverges. This issue is explored in Chapter 4 of this work.

The questions studied in the work are of scientific interest, are meaningful, have not been previously considered in the scientific literature and have scientific novelty.

## 2. GENERALIZATION WINTNER'S THEOREM

The Dirichlet convolution method is widely used among various methods for estimating the asymptotics of summation functions. This method is universal, as it can be used to estimate summation functions of additive, multiplicative and other arithmetic functions.

The idea of the method is as follows. Let there is a summation function $M(f, x) = \sum_{n \leq x} f(n)$ whose asymptotic estimate we want to obtain. Let's try to represent an arithmetic function in the form $f = f_0 * g = \sum_{d \mid n} g(d) f_0(n/d)$, where $f_0$ is a close to arithmetic function with known asymptotics $M_0(f, x) = \sum_{n \leq x} f_0(n)$ and $g$ a sufficiently small, in a certain sense, arithmetic function.

Then we obtain the following asymptotic behavior of the summation function $M(f, x) = \sum_{n \leq x} f(n)$ for $x \to \infty$:

$$M(f, x) = \sum_{n \leq x} f(n) = \sum_{n \leq x} \sum_{d \mid n} g(d) f_0(n/d) = \sum_{d \leq x} g(d) M_0(x/d) =$$

$$= \sum_{d \leq x} g(d)[M_0(x/d)] + O(\sum_{d \leq x} |g(d)|) = \sum_{n \leq x} g(n) M_o(x/n) + O(\sum_{n \leq x} |g(n)|). \qquad (2.1)$$

Assertion 1

Let the arithmetic function $f = f_0 * g$, where $f_0 = Id_s$ ($Id_s$ is the arithmetic function for which $Id_s = n^s, s = 0, 1, 2, ...k$) and $\sum_{n=1}^{\infty} \frac{|g(n)|}{n^{s+1}} < \infty$. Then the asymptotic behavior of the summation function $M(f, x) = \sum_{n \leq x} f(n)$ will be equal to:

$$M(f,x) = \sum_{n \leq x} f(n) = \frac{x^{s+1}}{s+1} \sum_{n=1}^{\infty} \frac{g(n)}{n^{s+1}} + o(x^{s+1}). \qquad (2.2)$$

Proof

Based on (1.6), we find the asymptotics $M_0(f,x) = \sum_{n \leq x} f_0(n)$ for $x \to \infty$:

$$M_0(n^s, x) = \sum_{n \leq x} n^s = \int_0^x t^s dt + O(x^s) = \frac{x^{s+1}}{s+1} + O(x^s). \qquad (2.3)$$

Let us use (2.1) and (2.3) to find the asymptotics of the summation function $M(f,x) = \sum_{n \leq x} f(n)$ for $x \to \infty$:

$$M(f,x) = \sum_{n \leq x} f(n) = \sum_{n \leq x} g(n) M_o(x/n) + O(\sum_{n \leq x} |g(n)|) = \frac{x^{s+1}}{s+1} \sum_{n \leq x} \frac{g(n)}{n^{s+1}} + O(|x^s \sum_{n \leq x} \frac{g(n)}{n^s}|) + O(\sum_{n \leq x} |g(n)|). \qquad (2.4)$$

Having in mind that the series $\sum_{n=1}^{\infty} \frac{|g(n)|}{n^{s+1}}$ converges, then the series $\sum_{n=1}^{\infty} \frac{g(n)}{n^{s+1}}$ converges, so when $x \to \infty$ based on (2.4) we obtain:

$$M(f,x) = \sum_{n \leq x} f(n) = \frac{x^{s+1}}{s+1} (\sum_{n=1}^{\infty} \frac{g(n)}{n^{s+1}} + O(|\sum_{n>x} \frac{g(n)}{n^{s+1}}|)) + O(|x^s \sum_{n \leq x} \frac{g(n)}{n^s}|) + O(\sum_{n \leq x} |g(n)|). \qquad (2.5)$$

On the other hand, since the series $\sum_{n=1}^{\infty} \frac{|g(n)|}{n^{s+1}}$ converges, therefore, based on the corollary from Kronecker's lemma (1.11), the following holds:

$$\sum_{n \leq x} |g(n)| = o(x^{s+1}). \qquad (2.6)$$

If $\sum_{n=1}^{\infty} \frac{g(n)}{n^{s+1}}$ - converges, then $\sum_{n>x} \frac{g(n)}{n^{s+1}} = O(1)$. In addition, in Chapter 4 it will be proven that if converges $\sum_{n=1}^{\infty} \frac{g(n)}{n^{s+1}}$ and $\sum_{n=1}^{\infty} \frac{g(n)}{n^s}$ - diverges, then based on (4.2) the following holds:

$$\sum_{n \leq x} \frac{g(n)}{n^s} = o(x). \qquad (2.7)$$

Based on (2.6) and (2.7), expression (2.5) for $x \to \infty$ can be written as:

$$M(f,x) = \sum_{n \leq x} f(n) = \frac{x^{s+1}}{s+1} \sum_{n=1}^{\infty} \frac{g(n)}{n^{s+1}} + o(x^{s+1}),$$

which corresponds to (2.2).

Let us note that at $s = 0$ expression (2.2) corresponds to the corollary from Wintner's theorem (1.5).

Let's look at an example for assertion 1.

Let it is necessary to find the asymptotic behavior of the summation function

$$M(\varphi, x) = \sum_{n \leq x} \varphi(n).$$

It is known that $\varphi = Id * \mu$, therefore $s = 1, g = \mu$. Let's check the condition $\sum_{n=1}^{\infty} \frac{|g(n)|}{n^{s+1}} = \sum_{n=1}^{\infty} \frac{|\mu(n)|}{n^2} = \frac{\zeta(2)}{\zeta(4)} < \infty$. Based on (2.2) then $x \to \infty$ we obtain:

$$M(\varphi, x) = \sum_{n \leq x} \varphi(n) = \frac{x^2}{2} \sum_{n=1}^{\infty} \frac{\mu(n)}{n^2} + o(x^2) = \frac{x^2}{2\zeta(2)} + o(x^2) = \frac{3x^2}{\pi^2} + o(x^2).$$

Another example for assertion 1.

Let it be necessary to find the asymptotic behavior of the summation function

$$M(\sigma_k, x) = \sum_{n \leq x} \sigma_k(n).$$

It is known that $\sigma_k = Id_k * 1$, therefore $s = k, g = 1$. Let's check the condition:

$$\sum_{n=1}^{\infty} \frac{|g(n)|}{n^{s+1}} = \sum_{n=1}^{\infty} \frac{1}{n^{k+1}} = \zeta(k+1) < \infty.$$

Based on (2.2) for $x \to \infty$ we obtain:

$$M(\sigma_k, x) = \sum_{n \leq x} \sigma_k(n) = \frac{x^{k+1}}{k+1} \sum_{n=1}^{\infty} \frac{1}{n^{k+1}} + o(x^{k+1}) = \frac{\zeta(k+1)x^{k+1}}{k+1} + o(x^{k+1}).$$

Now let's look at an example of linear asymptotics.

Let it be necessary to find the asymptotic behavior of the summation function

$$M(\frac{\varphi(n)}{n}, x) = \sum_{n \leq x} \frac{\varphi(n)}{n}.$$

It is known that $\frac{\varphi}{Id} = 1 * \frac{\mu}{Id}$, therefore $s = 0, g = \mu / Id$.

Let's check the condition: $\sum_{n=1}^{\infty} \frac{|g(n)|}{n} = \sum_{n=1}^{\infty} \frac{|\mu(n)|}{n^2} = \frac{\zeta(2)}{\zeta(4)} < \infty$.

Based on (2.2) for $x \to \infty$ we obtain:

$$M(\frac{\varphi(n)}{n}, x) = \sum_{n \leq x} \frac{\varphi(n)}{n} = x \sum_{n=1}^{\infty} \frac{\mu(n)}{n^2} + o(x) = \frac{x}{\zeta(2)} + o(x).$$

## 3. FINDING A CLASS OF ARITHMETIC FUNCTIONS THAT HAVE A LOGARITHMIC ASYMPTOTIC MEAN

Assertion 2

The summation function $\sum_{n \leq x} f(n)$ has a logarithmic asymptotic mean value if $f(n) = (\frac{A}{Id} * g)(n)$, where $A$ is a constant, and $g(n)$ satisfies the condition $\sum_{n=1}^{\infty} |g(n)| < \infty$.

Proof

We will use the Dirichlet convolution method and look for a function of the form $f = f_0 * g$. Let us take $f_0(n) = \frac{A}{n}$, therefore, based on (1.6):

$$M_0(f_0, x) = \sum_{n \leq x} f_0(n) = \sum_{n \leq x} \frac{A}{n} = A \ln x + A\gamma + O(1/n) = A \ln x + O(1), \qquad (3.1)$$

those. has a strong logarithmic asymptotic mean.

Then, based on (2.1), (3.1), we obtain the asymptotics $\sum_{n \leq x} f(n)$:

$$M(f, x) = \sum_{n \leq x} f(n) = \sum_{n \leq x} \frac{A}{n} * g(n) = \sum_{d \leq x} g(d) M_0(x/d) + O(\sum_{d \leq x} |g(d)|). \qquad (3.2)$$

Let's substitute (3.1) into (3.2) and get:

$$\sum_{n\leq x}f(n)=A\sum_{d\leq x}g(d)\ln(x/d)+O(\sum_{d\leq x}|g(d)|)=A\ln x\sum_{d\leq x}g(d)-A\sum_{d\leq x}g(d)\ln d+O(\sum_{d\leq x}|g(d)|). \quad (3.3)$$

Since the condition $\sum_{n=1}^{\infty}|g(n)|<\infty$ is satisfied, then $|g(n)|=O(1/n^{1+\xi})$ where $\xi>0$, therefore:

$$\sum_{n\leq x}|g(n)|\ln n\leq C\sum_{n\leq x}\frac{\ln n}{n^{1+\xi}}\leq C\sum_{n\leq x}\frac{\ln n}{n}=C\int_{2}^{x}\frac{\ln t\,dt}{t}+O(1)=C\ln\ln x+O(1). \quad (3.4)$$

Based on (3.4):

$$\sum_{n\leq x}|g(n)|\ln n=O(\ln\ln x). \quad (3.5)$$

Since $\sum_{n=1}^{\infty}|g(n)|<\infty$, then $\sum_{n=1}^{\infty}g(n)$ - converges and having in mind (3.3) and (3.5) we obtain:

$$\sum_{n\leq x}f(n)=A\ln x(\sum_{d=1}^{\infty}g(d)+o(1))+O(\ln\ln x)+O(1)=A\ln x\sum_{n=1}^{\infty}g(n)+o(\ln x) \quad (3.6)$$

Asymptotics (3.6) means that the summing function $\sum_{n\leq x}f(n)$ has a logarithmic asymptotic mean.

Let's look at an example for assertion 2.

Let $f_0(n)=\dfrac{3}{n}$ and $g(n)=\dfrac{\mu(n)}{n^2}$. Let's check the condition:

$$\sum_{n=1}^{\infty}|g(n)|=\sum_{n=1}^{\infty}\frac{|\mu(n)|}{n^2}<\infty.$$

Then, based on (3.6), we obtain the asymptotics of the summation function:

$$\sum_{n\leq x}\frac{3}{n}*\frac{\mu(n)}{n^2}=3\ln x\sum_{n=1}^{\infty}\frac{\mu(n)}{n^2}+o(\ln x)=\frac{3\ln x}{\zeta(2)}+o(\ln x).$$

Therefore, an arithmetic function $\dfrac{3}{n}*\dfrac{\mu(n)}{n^2}$ has a logarithmic asymptotic mean.

Let's consider another example for assertion 2.

Let $f_0(n) = \dfrac{1}{n}$ and $g(n) = \dfrac{\sigma(n)}{n^3}$. Let's check the condition:

$$\sum_{n=1}^{\infty} |g(n)| = \sum_{n=1}^{\infty} \frac{|\sigma(n)|}{n^3} < \infty ..$$

Then, based on (3.6), we obtain the asymptotics of the summation function:

$$\sum_{n \le x} \frac{1}{n} * \frac{\sigma(n)}{n^3} = \ln x \sum_{n=1}^{\infty} \frac{\sigma(n)}{n^3} + o(\ln x) = \zeta(2)\zeta(3) \ln x + o(\ln x).$$

Therefore, an arithmetic function $\dfrac{1}{n} * \dfrac{\sigma(n)}{n^3}$ has a logarithmic asymptotic mean.

## 4. ASYMPTOTICS OF SUMMATION FUNCTIONS WHEN CORRESPONDING DIRICHLET SERIES DIVERGES

Assertion 3

Let $f : \mathrm{N} \to C$ is an arithmetic function, $s$ is a complex number, the Dirichlet series $\sum_{n=1}^{\infty} \dfrac{f(n)}{n^s}$ converges, and the Dirichlet series $\sum_{n=1}^{\infty} \dfrac{f(n)}{n^{s-1}}$ diverges.

Then:

$$\sum_{n \le x} \frac{f(n)}{n} = o(x^{s-1}), \sum_{n \le x} \frac{f(n)}{n^2} = o(x^{s-2}), \dots, \sum_{n \le x} \frac{f(n)}{n^m} = o(x^{s-m}) \qquad (4.1)$$

at the value $x \to \infty$ and $\mathrm{Re}(s-k) > 0$, and

$$\sum_{n \le x} \frac{f(n)}{n^{s-k}} = o(x^k) \qquad (4.2)$$

at the value $x \to \infty$ and $\mathrm{Re}(s-k) > 0$.

Proof

We will use the discrete Abel transform to estimate the asymptotics of the summation function:

$$\sum_{n \le x} a_n f(n) = a(x) B(x) - \sum_{n \le x-1} (a_{n+1} - a_n) B(n). \qquad (4.3)$$

Let $B(x) = \sum_{n \leq x} f(n)$, $a_n = 1/n^k$, $k$ is a natural number, then based on (4.3) we obtain:

$$\sum_{n \leq x} \frac{f(n)}{n^k} = \frac{1}{x^k} \sum_{n \leq x} f(n) - \sum_{n \leq x-1} \left( \frac{1}{(n+1)^k} - \frac{1}{n^k} \right) \sum_{l \leq n} f(l). \tag{4.4}$$

Since $\sum_{n=1}^{\infty} \frac{f(n)}{n^s}$ - converges, and the Dirichlet series $\sum_{n=1}^{\infty} \frac{f(n)}{n^{s-1}}$ - diverges, based on the corollary to Kronecker's lemma, we obtain:

$$\sum_{n \leq x} f(n) = o(x^s), \sum_{l \leq n} f(l) = o(n^s). \tag{4.5}$$

Let's substitute (4.5) into (4.4) and get:

$$\sum_{n \leq x} \frac{f(n)}{n^k} = o(x^{s-k}) + \sum_{n \leq x-1} \frac{o(n^s) O(n^{k-1})}{n^k (n+1)^k} = o(x^{s-k}) + o\left( \sum_{n \leq x-1} n^{s-k-1} \right) = o(x^{s-k}). \tag{4.6}$$

Let us take the values $k = 1, 2, ..., m$ in (4.6) and obtain asymptotics (4.1) for $x \to \infty$ and $\mathrm{Re}(s-m) > 0$.

Now let's take $a_n = 1/n^{s-k}$, $k$ - a natural number, then based on (4.3) we get:

$$\sum_{n \leq x} \frac{f(n)}{n^{s-k}} = \frac{1}{x^{s-k}} \sum_{n \leq x} f(n) - \sum_{n \leq x-1} \left( \frac{1}{(n+1)^{s-k}} - \frac{1}{n^{s-k}} \right) \sum_{l \leq n} f(l). \tag{4.7}$$

Let us substitute (4.5) into (4.7) and obtain for $x \to \infty$ and $\mathrm{Re}(s-k) > 0$:

$$\sum_{n \leq x} \frac{f(n)}{n^{s-k}} = o(x^k) + \sum_{n \leq x-1} \frac{o(n^s) O(n^{s-k-1})}{n^{s-k}(n+1)^{s-k}} = o(x^k) + o\left( \sum_{n \leq x-1} n^{k-1} \right) = o(x^k),$$

which corresponds to (4.2).

Example for assertion 3.

It is known that the Dirichlet series $\sum_{n=1}^{\infty} \frac{\sigma_k(n)}{n^s} = \zeta(s)\zeta(s-k)$ at $\mathrm{Re}\, s > k+1$. Therefore, the series $\sum_{n=1}^{\infty} \frac{\sigma_k(n)}{n^s}$ converges at $s = k+2$ and diverges at $s = k+1$. Based on the corollary to Kronecker's lemma $\sum_{n \leq x} \sigma_k(n) = o(x^{k+2})$ in this case.

Using assertion 3 in this case:

$$\sum_{n\leq x}\frac{\sigma_k(n)}{n} = o(x^{k+1}),...,\sum_{n\leq x}\frac{\sigma_k(n)}{n^{k+1}} = o(x).$$

## 5. CONCLUSION AND SUGGESTIONS FOR FURTHER WORK

Next article will continue to study the asymptotic behavior of some arithmetic functions.

## 6. ACKNOWLEDGEMENTS

Thanks to everyone who has contributed to the discussion of this paper. The author is very pleased Professor Karras Meselem at Université de Khemis Miliana that showed interest in this article and read it in full before publication in the Archive.